\documentclass[final,hidelinks]{siamart190516}
\usepackage{amsmath,amssymb,mathtools,setspace,dsfont,hyperref,color,dsfont}
\newcommand{\N}{\mathbb{N}}
\newcommand{\Z}{\mathbb{Z}}
\newcommand{\R}{\mathbb{R}}

\newcommand{\Sb}{\mathbb{S}}

\newcommand{\E}{\mathbb{E}}
\newcommand\EN{\mathbb{E}\,\mathrm{N}}

\newcommand{\Rp}{\R_{+}}

\newcommand{\GL}{\mathrm{GL}}

\newcommand{\eps}{\varepsilon}
\newcommand{\s}{\sigma}

\renewcommand{\a}{\alpha}
\renewcommand{\b}{\beta}
\newcommand{\vol}{\mathrm{vol}}

\newcommand{\rk}{\mathrm{rank}}

\newcommand{\Prob}{\mathrm{Prob}}
\newcommand{\cone}{\mathrm{cone}}

\newcommand{\cM}{\mathcal{M}}

\DeclarePairedDelimiter\abs{\lvert}{\rvert}

\usepackage{color}
\definecolor{red}{rgb}{.7,0,0}

\usepackage{xcolor}
\renewcommand{\thispagestyle}[1]{} 

\newsiamremark{rem}{Remark}
\crefname{rem}{Remark}{Remarks}
\newsiamremark{qn}{Question}
\crefname{qn}{Question}{Questions}
\newsiamremark{conj}{Conjecture}
\crefname{conj}{Conjecture}{Conjectures}
\newsiamremark{ex}{Example}
\crefname{ex}{Example}{Examples}

\title{On the number of real zeros of random fewnomials\thanks{%Submitted to the editors DATE.
\funding{P.B. is supported by DFG grant BU 1371 2-2 and 
by the European Research Council (ERC) under the European's Horizon 2020 research and innovation programme 
(grant agreement No 787840). A.E. and J.T.-C. are supported by the Einstein Foundation, Berlin.}}}
\author{Peter B\"urgisser\thanks{Technische Universit\"at Berlin, Institut f\"ur Mathematik, Sekretariat MA 3-2, Strasse des 17. Juni 136, 10623, Berlin, Germany (\email{pbuerg@math.tu-berlin.de}, \email{ton-cue@math.tu-berlin.de})}
\and Alperen A. Erg\"ur\thanks{Carnegie Mellon University, School of Computer Science, 5000 Forbes Avenue, 15213 Pittsburgh, PA, USA. (\email{aergur@cs.cmu.edu})}
\and Josu\'e Tonelli-Cueto\footnotemark[2]
}

\begin{document}
\maketitle

\begin{abstract} 
Consider a system $f_1(x)=0,\ldots,f_n(x)=0$
of $n$ random real polynomial equations in $n$ variables,
where each $f_i$ has a prescribed set of exponent vectors described by a set
$A\subseteq \N^n$ of cardinality~$t$. 
Assuming that the coefficients of the $f_i$ are independent Gaussians of any variance, 
we prove that the expected number of zeros of the random system in the positive orthant is 
bounded from above by $\frac{1}{2^{n-1}}\binom{t}{n}$.
\end{abstract} 

%%%
\section{Introduction}

In many applications, 
we are faced with the problem of understanding or finding 
the (positive) real solutions of a system of multivariate polynomial equations, e.g., see
\cite{drton-sturmfels:09,horn-jackson:72,sottile-book:11}. 
Descartes' rule of signs~\cite[p.~42]{descartes}, dating from~1687,
is one of the oldest results providing information on this.
This rule directly implies that a real univariate polynomial with $t$ terms 
has at most $t-1$ positive zeros. 
In 1980, Khovanskii~\cite{kho:80} (see also \cite{kho:91}) 
obtained a far reaching generalization. 
He showed that a system 
$f_1(x)=0,\ldots,f_n(x)=0$ 
of $n$ real polynomial equations in~$n$ variables $x_1,\ldots,x_n$ 
involving $t$ distinct exponent vectors has no more than 
$$
 2^{\binom{t-1}{2}}\, (n+1)^{t-1}
$$
nondegenerate positive solutions.  
(A solution~$x$ is called nondegenerate if the derivative 
$D_xf$ is invertible. 
Khovanskii's result in fact allows for real exponents and he
even derived  more general quantitative bounds 
for broader families of analytic functions.) 
Following Khovanskii, one calls such a system a {\em fewnomial system} 
with $t$ exponent vectors. 
Like in Descartes' rule, the upper bound does not depend on the degrees, 
but only on the number of exponent vectors and the number of variables.
Khovanskii's bound was improved by Bihan and
Sottile~\cite{biha_so:07}, 
who proved the upper bound 
$\frac{e^2+3}{4}\, 2^{\binom{t-n-1}{2}}\, n^{t-n-1}$,
which is polynomial in $n$ 
when $k:= t-n$ is fixed;  
see also~\cite{sottile-book:11}. 

We note that these bounds are exponential in the number $t$ of exponent vectors. 
It is widely conjectured that these bounds are far from being optimal, 
but very little is known. The following question is a central open
problem in fewnomial theory \cite{rojas_phillipson}.

\begin{qn} \label{ahh}
Fix the number $n$ of variables. Is the number of nondegenerate positive solutions 
of a fewnomial system with $t$ exponent vectors bounded by a polynomial in $t$?
\end{qn}

This question is open even in the case of two variables
to the best of our knowledge; see~\cite{ko-po-ta:15}.
Several articles give an affirmative answer in special cases, 
such as the intersection of a line with the zero set of a bivariate $t$-nomial~\cite{avendano,paul}, 
or the intersection of two plane curves, where one is defined by a trinomial and 
the other by a $t$-nomial \cite[Cor.~16]{KPS:15}.
There is a very interesting connection to complexity theory \cite{koiran,bu-bri:18}. 

The current state of our understanding in fewnomial theory motivated us to investigate 
\cref{ahh} for random sparse polynomial systems.
Some of the motivation for this comes from Shub and Smale's well known work~\cite{Bez2}, 
in which for the first time a real probabilistic version of B\'ezout's theorem %over the reals 
was obtained.  
Our main result, \cref{main} below, implies that a natural, probabilistic version of 
\cref{ahh} has an affirmative answer for Gaussian random polynomials. 

In order to state this result, let us introduce some terminology. 
%We write $[n] := \{1,\ldots,n\}$. 
We denote the multiplicative group of positive real numbers by $\Rp$ and 
write $\Rp^n$ for the positive real orthant. 
Moreover, $x^\a:= x_1^{\a_1}\cdots x_n^{\a_n}$ stands for the monomial term  
with the exponent vector $\a=(\a_1,\ldots,\a_n)\in \Z^n$, where we abbreviate 
$x:=(x_1,\ldots,x_n)$. 
Fix a finite subset $A\subseteq\Z^n$ of cardinality~$t$ together with  
a map $\s\colon A \to \Rp$. We assign to this data the random polynomial system 
$$ 
 f_1(x) :=\sum_{\a \in A} \s(\a) \, \xi_{1,\a} x^{\alpha} ,\ldots,
 f_n(x) :=\sum_{\a \in A} \s(\a) \, \xi_{n,\a} x^{\alpha} ,
$$
where the $\xi_{i,\alpha} \sim \mathcal{N}(0,1)$ are independent identically distributed (i.i.d.) 
standard Gaussian random variables.
This amounts to considering random polynomial systems, 
%$f=(f_1,f_2,\ldots,f_n)$, 
where each $f_i$ has the support $A$, and the coefficients occurring in $f_i$ 
are independent centered Gaussian coefficients,  
whose variances are given by $\s(\alpha)^2$ for $\alpha\in A$. 
For a given {\em support} $A$ and a {\em system of variances} $\s$
we denote by $\EN(A,\s)$ the 
{\em expected number of nondegenerate zeros} in $\Rp^n$
of the random system $(f_1,\ldots,f_n)$. %$f(x)=0$ 

We can now state our main result. 

\begin{theorem}\label{main}
We have 
$$
 \EN(A,\s) \; \le\;  \frac{1}{2^{n-1}}\binom{t}{n} 
$$
for any support $A\subseteq\Z^n$ of cardinality~$t$ 
and any system of variances $\s\colon A \to \Rp$. 
\end{theorem}

\begin{rem}Our result holds in greater generality.
Consider the projection 
\[ \R^A\setminus \{0\} \longrightarrow \Sb(\R^A) \; , \;  y \mapsto y/\|y\| \]
to the unit sphere.  The proof only requires that for each~$i$, the image 
$\pi\big( (\xi_{i\a})_{\a \in A}\big)$ is uniformly distributed in the unit sphere of $\R^A$.
\end{rem}

\begin{rem}
By Markov's inequality, the probability that a random fewnomial system 
has a nondegenerate positive solution is bounded from above by $\EN(A,\s)$. 
In the situation, where $k=t-n$ is fixed and $n\to\infty$,  
\cref{main} implies that $\EN(A,\s) \to 0$ exponentially fast in~$n$.
Hence, in this situation, the system has no nondegenerate positive solution 
with overwhelming probability.

It may be interesting to compare the above with the deterministic polynomial upper bound 
in~\cite{biha_so:07} and with the following tight bounds. %observations. 
If $k=1$, there is at most one nondegenerate positive solution, cf.~\cite[Lemma~2]{itenberg-roy}.
Moreover, in the case $k=2$, there are at most $n+1$ nondegenerate positive solutions; 
see~\cite[Theorem~1.6]{rojas_phillipson} and~\cite{bihanbound}.
\end{rem}

To the best of our knowledge, the only previous results on random real fewnomial systems 
are by Malajovich and Rojas~\cite{malajovich,rojas_avg}.
These works provide upper bounds on the expected number of real zeros in terms of the (mixed) volume 
of the Newton polytopes of the $f_i$. Thus these bounds depend also on the degree, while our bound 
depends solely on the number of exponent vectors. 
We refer to Shiffman and Zelditch~\cite{shiff-zell:04,shiff-zell:11} for results on the distribution of the zeros 
of complex random fewnomials.  

There is a rich literature (see~\cite{bhar-reid:86})  
on the real zeros of random univariate polynomials, 
which mainly focuses on dense univariate polynomials, i.e., 
$n=1$ and $A=\{0,1,\ldots,d\}$. 
We only mention here Kac~\cite{kac:43}, who showed 
$\EN(A,\s) =  (\pi^{-1} + o(1)) \log d$ in the case $\s(\a)=1$, 
and Shub and Smale's paper~\cite{Bez2}, where 
$\EN(A,\s) = \frac12 d^{\frac{1}{2}}$ is proved when $\s(\a) = \binom{d}{a}^{\frac12}$. 
However, our focus is on the case of arbitrary supports $A$ 
(think of few exponent vectors of high degree), 
where apparently little is known. 
In the special case $n=1$, \cref{main} gives an upper of~$t$ which does not improve Descartes' rule. 
However, we can improve this result, 
%even though not for an arbitrary system of variances.
albeit for a specific system of variances only. 

\begin{theorem}\label{thm:uni}
In the univariate case, we have 
$$
 \EN(A,\mathds{1}) \; \le\;  \frac{2}{\pi}\sqrt{t}\log\,t
$$
for any support $A\subseteq\Z$ of cardinality~$t$
and the system of variances $\mathds{1}\colon \alpha\mapsto 1$.
\end{theorem}

The core idea of our approach can already be found in earlier work of Edelman and Kostlan~\cite{edel_kost}. 
Accordingly, one can express the expectation $\EN(A,\s)$ as the volume of a manifold, 
which can be seen as a real sparse analogue of the Veronese variety; 
see \cref{pro:Nequ}.  
While this proposition provides a beautiful geometric characterization of the expectation,  
the challenge consists of proving good estimations 
of the resulting integral. 
For instance, analytically proving any upper bounds on $\EN(A,\s)$---independent of the degree 
(which is possible thanks to Khovanskii's bound)---is a nontrival task. 
Our proof's main trick is a reduction %of the computation 
of the probabilistic upper bound to %the computation of 
a deterministic upper bound
for another explicit structured system.
This gives the proof a more combinatorial flavor and we leave it as a challenge to
find a more analytic proof. 
A natural next goal is to find an extension 
to the case where the fewnomials may have different supports.  
 
\bigskip

\noindent{\bf Acknowledgements.}
We are grateful to Felipe Cucker, Mario Kummer, and Gregorio Malajovich for 
helpful discussions. Special thanks go to Antonio Lerario for pointing out an error 
in a previous version of the paper. 
We are indebted to Irene\'{e} Briquel who contributed to \cref{thm:uni}.
We also thank Michael Joswig for useful discussions on discrete geometry, 
and J.~Maurice Rojas for passing his love of sparse polynomials on to his students. 
Finally, we thank the anonymous referees for detailed suggestions leading to an improvement 
of the presentation.

%%%
\section{Preliminaries}\label{se:prelim}

\subsection{Some integral geometry}

Random Gaussian vectors with independent coordinates are akin to uniform distributions over spheres, 
which brings us to the realm of integral geometry. 
Poincar\'e's formula for spheres (e.g., see \cite[Thm.~A.55]{cond}) implies the following result. 
(For a generalization to homogeneous spaces we refer to~\cite{howard:93} and \cite[Cor.~A.3]{bule:16}.) 
In what follows, $\vol_n$ denotes the volume measure induced by the usual Riemannian metric on~$\mathbb{S}^n$.

\begin{proposition}\label{cor:IG}
Let $\cM$ be an $n$-dimensional smooth submanifold of $\Sb^p$. Then we have
$$
 \E_{h_1,\cdots,h_n} \# (\cM \cap h_1 \Sb^{p-1} \cap\ldots\cap h_n \Sb^{p-1}) 
   = \frac{2\,\vol_n(\cM)}{\vol_n(\Sb^n)} ,
%   = \frac{\Gamma(\frac{n+1}{2})}{\pi^{\frac{n+1}{2}}} \, \vol_n(\cM) ,
$$
where the expectation is over independent, uniformly distributed 
$h_1,\ldots,h_n$ in the orthogonal group $O(p+1)$
defined with respect to the Haar measure.\hfill\qed 
\end{proposition}

\subsection{An integral formula for the expected number of zeros}

The core ideas of this subsection can be essentially found in~\cite{edel_kost}. 

Suppose we are given smooth and semialgebraic functions 
$\varphi_j\colon \Rp^n \to \R$, for $j=1,\ldots,t$ 
with $\varphi_1=\cdots=\varphi_t=0$ having no solutions in $\Rp^n$. 
%with empty zero set.
We consider random linear combinations 
\begin{equation}\label{eq:fxi}
 f_i(x) = \sum_{j=1}^t \xi_{ij} \varphi_j(x), \quad i=1,\ldots,n,
\end{equation}
where the $\xi_{ij} \sim \mathcal{N}(0,1)$ are i.i.d.\ standard
Gaussian random variables.
This defines a random smooth function $f\colon \Rp^n \to \R^n$, 
whose components are the $f_i$. 
Our goal is to study the number of nondegenerate zeros of $f$.
More specifically, 
we want to estimate the expectation $\E N$ of the random variable 
$$
 N(\xi) := \#\{ x \in \Rp^n \mid f(x) = 0,\  \det D_xf \ne 0 \} 
$$
taking values in $\N\cup\{\infty\}$. 
In order to do so, we consider the following smooth map 
$$
 \varphi \colon \Rp^n \to \R^t,\, x\mapsto (\varphi_1(x),\ldots,\varphi_t(x)) 
$$
and its scaled version 
$$
 \psi \colon \Rp^n \to \Sb^{t-1}, 
  \psi(x) := \frac{\varphi(x)}{\| \varphi(x)\|} ,
$$
which takes its values in the unit sphere. 

\begin{theorem}\label{pro:Nequ}
In the above setting, we have
$$
 \E N  = %\frac{n^{\frac{n}{2}}}{\pi^n}
\frac{2}{\vol_n(\Sb^n)} 
\, \int_{\Rp^n} \sqrt{\det\big((D_x\psi)^T D_x\psi\big)}\, dx .
$$
\end{theorem}
\begin{rem}
A tedious computation shows that this result in fact is~\cite[Thm.~3.3]{kostlan93}, 
see also~\cite[Theorem~7.1]{edel_kost}. 
Since these references contain only proof sketches, we provide a detailed proof.
%as in these references the arguments are only sketched.}
%in the proofs all details are missing.}
%~\cite{kostlan93} is missing all the details.
\end{rem}

\begin{proof}
We define the semialgebraic sets
\[
 V:=\{x\in\Rp^n\mid \rk D_x \psi < n\}\text{ and }U:=\{x\in\Rp^n\mid \rk D_x \psi = n\} .
\] 
Then we partition the open set $U$ into the semialgebraic sets
$$
 U_d := \{ x\in U \mid \#(\psi^{-1}(\psi(x))\cap U) = d \} ,
$$
for $d\in\N\cup\{\infty\}$. 
We associate with these semialgebraic sets the random variables
\[
N_V(\xi) := \#\{ x \in V \mid f(x) = 0,\  \det D_xf \ne 0 \} 
\]
and
$$
 N_d(\xi) := \#\{ x \in U_d \mid f(x) = 0,\  \det D_xf \ne 0 \} .
$$
Since $\{V,U_1,U_2,\ldots,U_\infty\}$ form a partition of $\Rp^n$, it suffices to prove that 
\begin{equation}\label{eq:gool1}
 \E N_V = %\frac{n^{\frac{n}{2}}}{\pi^n}
\frac{2}{\vol_n(\Sb^n)} \, \int_{V} \sqrt{\det\big((D_x\psi)^T D_x\psi\big)}\, dx ,
\end{equation}
and that for all $d\in\N\cup\{\infty\}$, 
\begin{equation}\label{eq:gool}
 \E N_d = %\frac{n^{\frac{n}{2}}}{\pi^n}
\frac{2}{\vol_n(\Sb^n)} \, \int_{U_d} \sqrt{\det\big((D_x\psi)^T D_x\psi\big)}\, dx .
\end{equation}

The right hand side of~\cref{eq:gool1} is zero since 
$\rk((D_x\psi)^T D_x\psi)=\rk D_x\psi<n$ for all $x\in V$. 
In order to prove~\cref{eq:gool1}, 
it is enough to show that $N_V=0$, which means that every zero $x\in V$ of the system $f$ 
is degenerate (i.e., $D_xf$ is singular). 
By \cref{eq:fxi}, we have $f(x) = [\xi_{ij}]\varphi(x)$ 
and so $D_xf=[\xi_{ij}]D_x\varphi$, 
for every $x\in\Rp^n$. 
By an explicit computation, we get 
\begin{equation}\label{eq:J}
 D_x\psi = \frac{1}{\|\varphi(x)\|} \Big( I - \psi(x)\psi(x)^T\Big) D_x \varphi .
\end{equation}
Suppose now that $x\in V$ satisfies $f(x)=0$. 
By~\cref{eq:J} and $\rk D_x\psi<n$,  
we either have $\rk D_x\varphi<n$, 
or there is some $v_x\in\mathbb{R}^n\setminus 0$ such that
$\varphi(x)=D_x\varphi\,v_x$. 
In the first case, $D_xf=[\xi_{ij}]D_x\varphi$ is singular. 
%hence $x$~is a degenerate zero of $f$. 
In the second case, $D_xf\,v_x=[\xi_{ij}]D_x\varphi\,v_x=[\xi_{ij}]\varphi(x)=f(x) =0$, 
hence $v_x\in\ker D_xf$ and $D_xf$ is singular as well.
We have thus shown that \cref{eq:gool1} holds.
%Hence, we have show that no $f$ can have nondegenerate zeroes inside $V$. 
%This means that $N_V=0$, proving~\cref{eq:gool1}.

For showing~\cref{eq:gool}, let $y_1,\ldots,y_t$ be new variables. 
We associate to the functions 
$f_i = \sum_{j=1}^t \xi_{i,j} \varphi_j(x)$ the linear forms 
$\ell_i := \sum_{j=1}^t \xi_{i,j} y_j$ and 
denote by $Z(\ell_1,\ldots,\ell_n)$ their zero set. 
So we have $f_i(x) = \ell_i(\varphi(x))$ for all~$x$. 
By the definition of~$U_d$, we have 
\begin{equation}\label{eq:NNN}
  \#\{ x \in U_d \mid f(x) = 0\} = d\, \# \big(\psi(U_d) \cap Z(\ell_1,\ldots,\ell_n) \big) .
\end{equation}

We first consider the case where $\dim\psi(U_d) =n$ and $d\in\N$. 
Using the stratification of semialgebraic sets into manifolds (cf.~\cite[Chap.~9]{bocr:98}), 
one shows that $\psi(U_d)$ contains a smooth $n$-dimensional submanifold $\cM_d$ of $\Sb^{t-1}$ 
such that $\dim(\psi(U_d)\setminus \cM_d) < n$.
By Sard's Theorem (cf.~\cite[\S A.2.4]{cond}),  
almost surely, the random hyperplanes $Z(\ell_1),\ldots,Z(\ell_n)$ 
intersect the $n$-dimensional manifold $\cM_d$ transversally, 
$$
 \psi(U_d) \cap Z(\ell_1,\ldots,\ell_n)  = 
  \cM_d \cap Z(\ell_1,\ldots,\ell_n) .
$$
Moreover, all the zeros of $f(x)=0$ in $U_d$ are nondegenerate.
With~\cref{eq:NNN} we conclude that, for almost all $\ell_1,\ldots,\ell_n$, 
$$
 N_d(\xi) = \#\{ x \in U_d \mid f(x) = 0 \} 
     = d\, \# \big(\psi(U_d) \cap Z(\ell_1,\ldots,\ell_n) \big) 
     = d\, \# \big(\cM_d \cap Z(\ell_1,\ldots,\ell_n) \big) .
$$
Therefore, applying \cref{cor:IG} 
to the manifold $\cM_d$, we obtain 
$$
 \E N_d = d\, \E \#\big(\cM_d \cap Z(\ell_1,\ldots,\ell_n)\big) = 
\frac{2d\, \vol_n(\cM_d)}{\vol_n(\Sb^n)} .
$$
Now we use that 
\begin{equation*}\label{eq:FL}
  \int_{U_d} \sqrt{\det\big((D_x\psi)^T D_x\psi\big)}\, dx 
   = \int_{y\in \cM_d} \# (\psi_d^{-1}(y) \cap U_d)   \, d\cM_d(y) 
   =  d\, \vol_n(\cM_d) ,
\end{equation*}
which follows from a slightly extended version of~\cite[Cor.~17.10]{cond}. 
(One can show that it does not matter that $\psi(U_d)$ 
may not be a manifold). 
This implies~\cref{eq:gool}.

In the case where $\dim\psi(U_d) <n$ and $d\in \N$, we write $\psi(U_d)$ as a
union of smooth manifolds of dimension less than~$n$. Then, using Sard's 
Theorem, we see that~\cref{eq:gool} trivially holds in the form $0=0$. 

To complete the proof, it suffices to show $U_\infty$ is empty. 
By way of contradiction, assume $x\in U_\infty$. 
By the Constant Rank Theorem~\cite[Theorem~4.12]{LEE}, 
the fiber $\psi^{-1}(\psi(x))\cap U$ is a zero-dimensional subset of the open set $U$, 
since $\rk D_p\psi=n$ for all $p\in U$. However, $\psi^{-1}(\psi(x))\cap U$ is semialgebraic and hence finite, 
since any zero-dimensional semialgebraic set is finite; cf.~\cite{bocr:98}. 
This contradicts $x\in U_\infty$, completing the proof.
\end{proof}

\section{Special systems with few terms}\label{se:special-systems}

We prove a deterministic result (\cref{pro:special-system}) on the real zeros 
of a particular class of sparse systems that involve 
square roots and thus goes slightly beyond our polynomial 
setting. This result will allow us to prove the following inequality of integrals 
that will be instrumental in the proof of \cref{main}.

\begin{proposition}\label{cor:alp-trick}
Consider the function $f\colon\Rp^n\to\Rp$ defined by 
$f(x) := \left(\sum_{j=1}^m c_j^2 x^{2\b_j}\right)^{\frac12}$, 
where $\b_1\ldots,\b_m \in\Z^n$ and $c_1,\ldots,c_m \in\R^*$.
Take $\a_1,\ldots,\a_n\in\Z^n$ and $\s_1,\ldots,\s_n>0$, 
and define the functions 
$\varphi\colon\Rp^n \to \Rp^{n+1}$ 
and $\psi\colon\Rp^n \to S^n$ by 
$$
 \varphi(x) := (\s_1 x^{\a_1}, \ldots, \s_n x^{\a_n}, f(x)) ,\quad 
 \psi(x) : = \frac{\varphi(x)}{\|\varphi(x)\|} .
$$
Then we have 
$$
\frac{1}{\vol_n(\Sb^n)} \,
  \int_{\Rp^n}  \sqrt{\det\big((D_x\psi)^T D_x\psi\big)}\, dx \ \le\ \frac{1}{2^{n}}.
$$
\end{proposition}

For proving \cref{cor:alp-trick} we need the following lemma.

\begin{lemma}\label{pro:special-system}
%Let $f(x)$ be the function from \cref{pro:special-system}, but with $\b_1\ldots,\b_m \in\R^n$.  
Following the notation of \cref{pro:special-system}, 
we allow any $\b_1\ldots,\b_m \in\R^n$ in the definition of $f$. 
Moreover, let $\a_1,\ldots,\a_n\in\R^n$. 
Then, for any $[\lambda_{ij}] \in \GL_n(\R)$, 
the system 
$$
 \sum_{j=1}^n \lambda_{ij} x^{\a_j} = f(x) ,\quad i=1,\ldots,n, 
$$
has at most two nondegenerate zeros in $\Rp^n$.
\end{lemma}

We recall the following fact about changing variables that we will use in the proof of the lemma.
Suppose $a_1,\ldots,a_n$ is a basis of~$\R^n$. 
Then 
\begin{equation}\label{eq:trafoexp}
 \Rp^n \to \Rp^n,\, x\mapsto (x^{a_1},\ldots,x^{a_n}) 
\end{equation}
is a diffeomorphism.
Indeed, via the group isomorphism $\R^n\to\Rp^n,\, y\mapsto \exp(y)$ and its inverse 
$\Rp^n\to\R^n,\, x\mapsto (\log x_1,\ldots,\log x_n)$,
this turns \cref{eq:trafoexp} into the linear isomorphism 
$\R^n\to\R^n,\, y\mapsto (a_1^T y,\ldots,a_n^Ty)$.

\begin{proof}[Proof of \cref{pro:special-system}]
We divide into cases, depending on the rank~$k$ of the linear span of $\alpha_1,\ldots,\alpha_n$.

In the case $k=n$, using the transformation as in~\cref{eq:trafoexp}, 
we can assume without loss of generality that 
$\a_1,\ldots,\a_n$ is the standard basis of $\R^n$; 
thus we study the positive zeros of a system 
\begin{equation}\label{eq:uno}
 \sum_{j=1}^n \lambda_{ij} x_j = f(x) ,\quad i=1,\ldots,n .
\end{equation}
Subtracting the $n$th equation from the others gives the system
$\sum_{j=1}^n (\lambda_{ij} -\lambda_{nj}) x_j = 0$, $i=1,\ldots,n-1$,
which has a one dimensional solution space $\R\xi$, for some nonzero $\xi\in\R^n$.  
We can assume that $\xi\in\Rp^n$ since otherwise the system \cref{eq:uno} has no solution in $\Rp^n$. 

Plugging in $x=s\xi$ with unknown $s\in\Rp$ into \cref{eq:uno}, 
we obtain by squaring the first equation
$$
 s^2 \left(\sum_{j=1}^m \lambda_{1j} \xi_{j}\right)^2 - \sum_{j=1}^m c_j^2 \xi^{2\b_j} s^{2 g_j} = 0,
$$
where $g_j$ denotes the sum of the components of $\beta_{j}$. 
We apply now Descartes' rule to this univariate polynomial in $s$ (with possibly real exponents). 
Since the sequence of coefficients has at most two sign changes, there are at most 
two positive real zeros, provided the polynomial does not vanish altogether. 
In the latter case, all the $g_j$ equal 1 and we have $f(s\xi) = sf(\xi)$.
The system~\cref{eq:uno} then becomes the system $s\sum_j\lambda_{ij}\xi_j=sf(\xi)$ in $s$, for $i=1,\ldots,n$,
whose solution set either is empty or all of $\Rp$.
But then all solutions of the original system~\cref{eq:uno} are degenerate.

In the case $k<n$, using the transformation as in~\cref{eq:trafoexp}, 
we can assume without loss of generality that 
$\a_1,\ldots,\a_n$ lies in $\R^k\times 0^{n-k}$.
%the linear span of the first standard basis vectors.
Subtracting the $n$th equation of the original system from the others gives the system
$\sum_{j=1}^n (\lambda_{ij} -\lambda_{nj}) x^{\alpha_j} = 0$, $i=1,\ldots,n-1$. 
Since $x^{\a_j}$ only depends on $x_1,\ldots,x_k$, this is a system of $n-1$ equations 
in $k$ variables, having $n$ exponent vectors.
If $k<n-1$, this system only has degenerate solutions and 
we are done.
If $k=n-1$, we are faced with a system of $n-1$ polynomials in $n-1$ variables having $n$ exponent vectors. 
It is known that such systems have at most one nondegenerate positive real solution~$\xi$; 
cf.~\cite[\S3.2.1]{sottile-book:11}.
Substituting 
$(\xi_1,\ldots,\xi_{n-1},x_n)$ with unknown~$x_n$ in the first equation,  
and using $\a_{j,n}=0$, we obtain after squaring 
$$
 \left(\sum_{j=1}^{n} \lambda_{1j} \xi_{1}^{\a_{j,1}} \cdots \xi_{n-1}^{\a_{j,n-1}}\right)^2  
 - \sum_{j=1}^m c_j^2 \xi_1^{2\b_{i,1}} \cdots \xi_{n-1}^{2\b_{i,n-11}}  x_n^{2 \b_{j,n}} = 0.
$$
By Descartes' rule, this polynomial in~$x_n$ (with possibly real exponents) has at most 
one positive zero,
%two positive zeros, 
unless it vanishes altogether, 
in which case all solutions of the original system are degenerate.

Summarizing, we have shown that in all cases, the system has at most two nondegenerate solutions. 
\end{proof}

\begin{rem}
\cref{pro:special-system} is optimal: 
if we take $f(x):= (1+ \frac15 x_1^2x_2^2)^{\frac12}$, then the system 
$x_1=f(x), x_2=f(x)$ has two positive real solutions. 
\end{rem}

\begin{proof}[Proof of \cref{cor:alp-trick}]
We begin with the following general observation: let 
$v_1,\ldots,v_n,w$ be independent standard Gaussian vectors in $\R^n$. 
Then 
$$
 \sum_{\eps\in\{-1,1\}^n} \Prob \big\{ w \in \cone(\eps_1 v_1,\ldots, \eps_n v_n) \big\}
 =\Prob\big\{ w \in \mathrm{span}(v_1,\ldots, v_n) \big\} = 1 .
$$
By symmetry, the probabilities do not depend on $\eps$, so that we get 
\begin{equation}\label{eq:22n}
 \Prob\{w \in \cone(v_1,\ldots,v_n) \} = 2^{-n} .
\end{equation}

Suppose now~$\xi_{ij}$ are independent standard Gaussian random variables. 
If the random system 
\begin{equation}\label{eq:systemspecial}
  \xi_{i0}f(x) + \sum_{j=1}^n \xi_{ij} \s_j x^{\a_j} =0, \quad i=1,\ldots,n,
\end{equation}
has a positive solution, then $-\xi_{0}$ is a positive linear combination of $\xi_1,\ldots,\xi_n$ and so
\[
 -\xi_{0}\in \cone(\xi_1,\ldots,\xi_n).
\]
The probability of the later is at most $2^{-n}$ by \cref{eq:22n}. 
Hence the probability that the random system \cref{eq:systemspecial} has a positive root is also bounded 
from above by $2^{-n}$. 

By \cref{pro:special-system}, the maximum number of positive nondegenerate zeros is at most two. 
Therefore, the expected number of nondegenerate solutions 
of \cref{eq:systemspecial} is bounded from above by $2 \cdot 2^{-n}$.
The assertion follows now by \cref{pro:Nequ}.
\end{proof}

%%%
\section{Proof of \texorpdfstring{\cref{main}}{Theorem~1.2}}\label{se:proofs}

We fix a support $A\subseteq\Z^n$ of cardinality~$t$ and $\s\colon A\to\Rp$.
Consider the following map of Veronese type: 
\begin{equation}\label{eq:veronese}
 v_{A,\s}\colon \Rp^n \to \Rp^A,\, x\mapsto (\s(\a)x^\a)_{\a\in A} ,
\end{equation}
and its scaled version
$$
 \gamma_{A,\s}\colon \Rp^n \to \Sb(\R^A),\, 
  \gamma_{A,\s}(x) := \frac{v_{A,\s}(x)}{\|v_{A,\s}(x)\|} ,
$$
which takes its values in the unit sphere. 

As all the functions are semialgebraic, 
we can apply \cref{pro:Nequ}, which implies that 
$$
 \EN(A,\s) = \frac{2}{\vol_n(\Sb^n)} \, \int_{\Rp^n}  
  \sqrt{\det ( (D_x\gamma_{A,\s})^T D_x\gamma_{A,\s})}\, dx .
$$
We abbreviate 
$M(x):= D_x \gamma_{A,\s}$.
For $I\subseteq A$ with $|I|=n$ 
we denote by $M_I(x)$ the square submatrix of $M(x)$ 
obtained by selecting the rows with index in $I$. 
The Cauchy-Binet formula~\cite[Theorem~2.3]{CBref}, combined with the elementary inequality $\|u\|_2 \le \|u\|_1$, gives 
\begin{equation} \label{CB}
\sqrt{\det(M(x)^TM(x))} = \left(\sum_{\abs{I}=n} (\det M_I(x))^2 \right)^{\frac12} 
  \ \le\ \sum_{\abs{I}=n} |\det M_I(x) | .
\end{equation}
Therefore,
$$
 \EN(A,\s) \; \le\; \sum_{\abs{I}=n}
\frac{2}{\vol_n(\Sb^n)} \, \int_{\Rp^n}  |\det M_I(x) |\, dx .
$$
It suffices to prove that 
\begin{equation}\label{eq:SInt}
\frac{1}{\vol_n(\Sb^n)} \,
 \int_{\Rp^n} |\det M_I(x)|\, dx \ \le\ \frac{1}{2^{n}} , 
\end{equation}
since there are $\binom{t}{n}$ summands. 

For showing this, we put 
$I=\{\a_1,\ldots,\a_n\}$ and $\s_i :=\s(\a_i)$. 
We then apply \cref{cor:alp-trick}
to the function
$\varphi\colon\Rp^n \to \Rp^{n+1}$ defined by 
$$
 \varphi(x) := (\s_1 x^{\a_1}, \ldots, \s_n x^{\a_n}, f(x)) ,
$$
where 
$$
 f(x) := \left(\sum_{\a\in A \setminus I} \s(\a)^2 x^{2\a}\right)^\frac12 .
$$
Note that 
$$
 \|\varphi(x)\|^2 = \sum_{i=1}^n \s_i^2 x^{2\a_i} + \sum_{\a\in A\setminus I} \s(\a)^2 x^{2\a} 
   = \sum_{\a \in A} \s(\a)^2 x^{2\a} = \|v_{A,\s}(x)\|^2 .
$$
Moreover, the $i$th coordinate of the scaled function $\psi(x) := \varphi(x)/\|\varphi(x)\|$ 
satisfies for $1\le i \le n$, 
$$
 \psi_i(x) = \frac{\varphi_i(x)}{\|\varphi(x)\|} = \frac{\s_i x^{\a_i}}{\|\nu_{A,\s}(x)\|} , 
$$
which is the $\a_i$th coordinate of $\gamma_{A,\s}(x)$.   
Therefore, 
$M_I(x) = [ \partial_{x_j} \psi_i ]_{i,j \le n}$.  
The Cauchy-Binet formula implies that 
$$
 |\det M_I(x)| \; \le \;  
\sqrt{\det\big((D_x\psi)^T D_x\psi\big)} .
$$ 
\cref{cor:alp-trick} gives 
$$
\frac{1}{\vol_n(\Sb^n)} 
 \, \int_{\Rp^n}  \sqrt{\det\big((D_x\psi)^T D_x\psi\big)}\, dx \ \le\ \frac{1}{2^{n}} .
$$
Combining the above shows~\cref{eq:SInt} and completes the proof 
of \cref{main}. \qed

\section{Proof of \texorpdfstring{\cref{thm:uni}}{Theorem~1.5}}\label{sec:uni}

Given a random polynomial $f$ with support $A$, $g:=f(1/x)$  is a random Laurent polynomial with support $-A$,  
whose expected number of zeros in $(0,1)$ is precisely the expected number of zeros of $f$ in $(1,\infty)$. 
Therefore, it is enough to bound the expected number of zeros $\EN_{(0,1)}(A,\mathds{1})$ 
in the interval $(0,1)$ by $\frac{1}{\pi}t^{1/2}\log\,t$ 
for a random polynomial with arbitrary support $A$ of size $t$. 
Moreover, since multiplying by $x^k$ does not alter the number of zeros in $(0,1)$ of a polynomial, 
we can assume without loss of generality that $0\in A\subseteq \N$.

We observe that \cref{pro:Nequ} holds for any open subset of $\R_+^n$ with the same proof.
Hence
\[\EN_{(0,1)}(A,\mathds{1})=\frac{1}{\pi} 
 \, \int_0^1  \|\psi'(x)\|\, dx\]
where $\varphi(x)=(x^\alpha)_{\alpha\in A}$ and $\psi:=\varphi/\|\varphi\|$. 
%The simplification of the integrand happens because $\psi'$ is a vector in the univariate case. 
By~\cref{eq:J}, $\psi'(x)=\|\varphi(x)\|^{-1}P_x\varphi'(x)$ where $P_x$ is the orthogonal projection onto the orthogonal complement of $\psi(x)$, and so
$$
\|\psi'(x)\| =\frac{\|P_x\varphi'(x)\|}{\|\varphi(x)\|}\leq  \frac{\|\varphi'(x)\|}{\|\varphi(x)\|} .
$$
Hence,
\begin{equation}\label{eq:inequni}
\|\psi'(x)\|\leq \sqrt{t}\,\frac{\|\varphi'(x)\|_1}{\|\varphi(x)\|_1} 
 = \sqrt{t}\,\left(\ln\|\varphi(x)\|_1\right)' ,
\end{equation}
using the standard inequalities 
between the $\ell_1$ and $\ell_2$ norms.
Finally, we obtain by integrating 
$$
  \int_0^1\|\psi'(x)\|\, dx \le 
   \sqrt{t} \, \left( \ln \|\phi(1)\|_1 - \ln \|\phi(0)\|_1 \right)   
  \ \le\ \sqrt{t} \, \ln t,
$$
since $0\in A\subseteq \N$, which gives the desired result.
\hfill\qed

\bibliographystyle{siamplain}
\bibliography{fewnomialrefs}

\end{document}